\documentclass[12pt,a4paper]{article}
\usepackage{amssymb,latexsym}
\usepackage{authblk}

\newcommand{\carac}[1]{{\bf 1}_{ \{#1\}}}
\newcommand{\mulength}{\nu_*}
\newcommand{\muemetteur}{\nu}
\newcommand{\mucompteurfix}{m}
\newcommand{\muunif}{m}
\newcommand{\muoverflow}{\eta}
\newcommand{\muglob}{\mu^H}
\newcommand{\muunf}{\mu^F_l}
\newcommand{\mudeuxf}{\mu^{F}_l}
\newcommand{\muf}{\mu^{F}}

\newcommand{\ZZ}{\mathbb{Z}}
\newcommand{\NN}{\mathbb{N}}

\newcommand{\EE}{\mathbb{E}}

\newtheorem{thm}{Theorem}

\newtheorem{lem}{Lemma}
\newtheorem{conj}{Conjecture}
\newtheorem{cor}{Corollary}
\newtheorem{pro}{Proposition}
\newtheorem{remark}{Remark}

\def\OE{\overline{E}} 

\begin{document}


\title{On a zero speed sensitive cellular automaton \footnote{Published in Nonlinearity 20 (2007) 1–19 }}

\author[1]{Xavier Bressaud }
\author[2]{Pierre Tisseur }
\affil[1]{ Université Paul Sabatier
Institut de Math\'ematiques de Toulouse
, France  \thanks{E-mai adress: \texttt{bressaud@math.univ-toulouse.fr }}}
\affil[2]{ Centro de Matematica, Computa\c{c}\~ao e Cogni\c{c}\~ao, Universidade Federal do ABC, Santo Andr\'e, S\~ao Paulo, Brasil \thanks{E-mai adress: \texttt{pierre.tisseur@ufabc.edu.br}}}
 
\sloppy 
\date{12 october 2006 }

\maketitle

\begin{abstract}
Using an unusual, yet natural invariant measure we show that there exists a
sensitive cellular automaton whose perturbations propagate at
asymptotically null speed for almost all configurations.
More specifically, we prove that Lyapunov Exponents measuring
pointwise or average
linear speeds of the faster perturbations are equal to zero. We show
that this implies the nullity of the measurable
entropy.
The measure $\mu$ we consider gives the $\mu$-expansiveness property
to the automaton. It is constructed with respect to a factor
dynamical system based on simple ``counter dynamics''.
As a counterpart,  we prove that in  the case of positively expansive automata,
the perturbations move at positive linear speed over all the configurations.
\end{abstract}


  \section{Introduction}

 A one-dimensional cellular automaton is a discrete mathematical
 idealization of a space-time
 physical system. The space $A^\ZZ$ on which it acts is the set of
 doubly infinite sequences of elements of a finite set $A$; it is called
 the {\em configuration space}. The discrete time is
 represented by the action of a cellular automaton $F$ on this space.
 Cellular automata are a class of  dynamical systems on which
 two different kinds of measurable entropy can be considered: the entropy with respect
 to the shift $\sigma$ (which we call spatial) and the entropy with respect to $F$
 (which we call temporal).
 The temporal entropy depends on the way the automaton "moves" the
 spatial entropy using a local
 rule on each site of the configuration space.
 The propagation speed of the different one-sided configurations, also called perturbations in this case,
 can be defined on a specific infinite
 configuration, or as an average value on the configuration space
 endowed with a probability measure.
 We can consider perturbations moving from the left to the right or from the right to the left side of
 the two sided sequences.
 Here we prove that the perturbations (going to the left or to the right) move at a positive speed on all
 the configurations for a
 positively expansive cellular automata and that in the sensitive case,
 there exist automata
 with the property that for almost all the configurations, the
 perturbations can move to infinity but at
 asymptotically null speed.

 Cellular automata can be  roughly divided into two classes: the class
 of automata which have
 {\em equicontinuous points} and the class of {\em sensitive} cellular
 automata (K\r{u}rka introduces a more precise classification in
 \cite{Ku94}). This partition of CA into ordered ones and disordered ones
 also
 corresponds to the cases where the perturbations cannot move to infinity
 (equicontinuous class) and to the cases
 where there always exist perturbations that propagate to infinity.

 The existence of equicontinuity points is equivalent to the existence
 of a so-called "blocking
 word" (see \cite{Ku94}, \cite{ti2000}). Such a word stops the
 propagation of perturbations so roughly speaking  in the non sensitive case,
 the speed of propagation is equal to zero for all the points which
 contains infinitely many occurrences of "blocking words".
 For a sensitive automaton
 there is no such word, so that perturbations may go to infinity.
 However, there are few results about the speed of propagation
 of perturbations in sensitive cellular automata except for the
  positively expansive subclass.

In \cite{Sh91}, Shereshevsky gave a first formal definition of these
speeds of propagation (for each point and for almost all points) and called them Lyapunov
 exponents because of the  analogy (using an appropriate metric) with the well known
 exponents of the differentiable dynamical systems.

 Here we use a second definition of these discrete Lyapunov exponents
 given by Tisseur \cite{ti2000}.
 The two definitions are quite similar, but for each point, the Shereshevsky's one  uses a
 maximum value on the
 shift orbit. For this reason the Shereshevsky's exponents (pointwise or global) can not see the "blocking
 words" of the equicontinuous points and could give a  positive value
 (for almost all points)
to the speed on
 these points.
 Furthemore, in our sensitive example (see Section 5), for almost all  infinite
 configurations,
there always exists some  increasing (in size) sequences of finite configurations
where the speed of propagation is not asymptoticaly null which implies that the initial definition
gives a positive value to the speed and does not take into account a part of the  dynamic of the
measurable dynamical systems.

Using the initial definition
 of Lyapunov exponents due to Shereshevsky \cite{Sh91}, Finelli,
 Manzini, Margara (\cite{FMM98})
 have shown that positive expansiveness implies positivity
 of the Shereshevski pointwise Lyapunov exponents at all points.


 Here we show that the statement of Finelli, Manzini, Margara still
 holds for our definition of pointwise
 exponents and the main difference between the two results is that we obtain that the exponents
are positive for all points using a liminf rather that a limsup.

 \begin{pro}
 \label{theorem1}
 For a positively expansive cellular automaton $F$ acting on $A^\ZZ$,
 there is a constant $\Lambda >0$ such that, for all $x \in X$,
 $$\lambda^+(x)\ge \Lambda \hbox{ and } \lambda^-(x)\ge \Lambda , $$
 where $\lambda^+(x)$ and $\lambda^-(x)$ are respectively the right and left pointwise Lyapunov exponents.
 \end{pro}
 The first part of the proof uses standard compactness arguments. The
 result is stronger and the
 proof is completely different from the one in \cite{FMM98}. This result
 is called
 Proposition~\ref{positiveexpansive} in
Section~\ref{sectionpositivelyexpansive} and it is stated for all $F$-invariant subshifts $X$.

 \paragraph{}
 Our main result concerns sensitive automata and average Lyapunov exponents ($I_\mu^+$ and
 $I_\mu^-$).
 We construct a sensitive cellular automaton $F$ and a
 $(\sigma,F)$-invariant measure $\mu^F$ such that
 the average Lyapunov exponents $I_{\mu^F}^\pm$ are equal to zero.

 By showing (see Proposition 3) that the nullity of the average Lyapunov exponents implies
 that the measurable entropy is equal to zero, we obtain that our particular
 automaton have null measurable entropy $h_{\mu^F}(F)=0$.

 We also prove that this  automaton
 is not only sensitive but $\mu^F$-expansive which is a measurable equivalent to
 positive expansiveness introduced by Gilman in \cite{Gil87}.
 So even if this automaton is very close to positive expansiveness (in the
 measurable sense), its
 pointwise Lyapunov exponents are equal to zero almost everywhere (using
  Fatou's lemma) for a "natural" measure $\mu^F$ with positive entropy under the
 shift. The  $\mu^F$-expansiveness means
 that "almost all perturbations" move to infinity and the Lyapunov
 exponents represent the speed of the
 faster perturbation, so in our example almost all perturbations move to
 infinity at
 asymptoticaly null speed.

 In view of this example,  Lyapunov exponents or average speed of
 perturbations appear
 a useful tool for proving
 that a cellular automata has zero measure-theoretic entropy.

 The next statement gathers the conclusions of Proposition \ref{i=0},
 Proposition \ref{sensitivity}, Lemma \ref{proprietesmesures} ,
 Proposition \ref{shiftEntr}, Proposition \ref{ii=0}, Corollary 1, Remark 6.
 \begin{thm}\label{thm2}
 There exists a sensitive cellular automaton $F$ with the following
 properties:
 there exists a $(\sigma,F)$-invariant measure $\mu^F$ such that $h_{\mu^F}
 (\sigma )>0$ and the
 Lyapunov exponents are equal to zero, i.e., $I^\pm_{\mu^F} = 0$,
 which implies that
 $h_{\mu^F} (F)=0$. Furthemore this automaton $F$ has the
 $\mu^F$-expansiveness property.
 \end{thm}

 Let us describe the dynamics of the cellular automaton $F$ and the
 related
 "natural" invariant measure $\mu^F$ that we consider.

 In order to have a perturbation moving to infinity but at a
 sublinear speed, we define a cellular
 automaton with an underlying  "counters dynamics", i.e. with a factor
 dynamical system based on "counters dynamics".

 Consider this factor dynamical system of $F$ and
 call a "counter" of size $L$ a  set $\{ 0, \ldots, L-1\} \subset
 \NN$. A
 trivial dynamic on this finite set is addition of $1$ modulo $L$.
 Consider a bi-infinite sequence of counters (indexed by $\ZZ$). The
 sizes of the counters will be chosen randomly and unboundly. If at each
 time step every counter is increased by one, all counters count at
 their own rythm,
 independently. Now we introduce a (left to right) interaction between
 them.
 Assume that each time a counter reaches the top (or passes through
 $0$; or is at $0$),
 it gives an overflow to its right neighbour. That is, at each time, a
 counter increases
 by $1$ unless its left neighbour reaches the top, in which case it
 increases by $2$.
 This object is not a cellular automaton because its state space is
 unbounded. However, this rough definition should be enough to suggest
 the idea.

 $\bullet$ The dynamics is sensitive. Choose a configuration, if we
 change
 the counters to the left of some coordinate $-s$, we change the
 frequency of apparition of
 overflows in the counter at position $-s$ (for example this happens if
 we put larger and larger counters).
 Then the perturbation eventually appears at the coordinate $0$.

 $\bullet$ The speed at which this perturbation propagates is controlled
 by the sizes of the
 counters. The time it takes to get through a counter is more or less
 proportional to the size of this counter (more precisely, of the
 remaining time before it reaches $0$ without an overflow). So a good
 choice of the law of the
 sizes of the counters allows us to control this mean time. More
 specifically, we prove that
 that with high probability, information will move slowly. If the size
 of the counters were bounded the speed would remain linear.

 In the cellular automaton $F$,
 we "put up the counters" horizontally: we replace a counter of size
 $L = 2^l$ by a
 sequence of $l$ digits and we separate sequences of digits by a special
 symbol, say $E$. Between two $E$s the dynamics of a counter is
 replaced by an odometer with overflow transmission to the right. More
 precisely,
 at each step the leftmost digit is increased by $1$  and
 overflow is transmitted.

 Note that to model the action of a cellular automaton we
 need to introduce in the factor dynamics a countdown which starts when
 the counter reaches the top. The end of the countdown corresponds to
  the transmission of the  overflow.
 When the countdown is running,
 the time remaining before the emission of the overflow does not
 depend on a possible overflow emmited by a neighbouring counter.
 Nevertheless the effect of this overflow will affect the start of the
 next countdown.

 Finally, we construct an invariant measure based on Cesaro means of the sequence $(\mu\circ F^{n})$
 where $\mu$ is a measure defined thanks to the counter dynamic of the factor dynamical system .


\section{Definitions and notations}

\subsection{Symbolic systems and  cellular automata}

Let $A$ be a finite set or alphabet.
Denote by $A^*$ the set of all concatenations of letters in $A$.
$A^\ZZ$ is the set of bi-infinite sequences
$x=(x_i)_{i\in\ZZ}$ also called \emph{configuration space}.
For $i\le j$ in $\ZZ$, we denote by $x(i,j)$ the word $x_i\ldots x_j$ and
by $x(p,\infty )$ the infinite sequence $(v_i)_{i\in\NN}$ with $v_i=x_{p+i}$.
For $t\in\NN$ and a word $u$ we
call \emph{cylinder} the set $[u]_t=\{x\in A^\ZZ : x(t,t+|u|)=u\}$.
%
The configuration set $A^\ZZ$ endowed with the product topology
is a compact metric space.
A metric compatible with this topology can be defined by the distance
$d(x,y)=2^{-i}$,
where $i=\min\{|j|:\ x(j)\ne y(j)\}$.

The shift $\sigma \colon A^\ZZ\to A^\ZZ$ is defined by
$\sigma (x_i)_{i\in \ZZ}=(x_{i+1})_{i\in \ZZ}$.
The dynamical system $(A^\ZZ ,\sigma )$ is called the \emph{full shift}.
A \emph{subshift} $X$ is a closed shift-invariant
subset $X$ of $A^\ZZ$ endowed with the shift
$\sigma$. It is possible to identify $(X,\sigma )$ with the set $X$.

Consider a probability measure $\mu$ on the Borel sigma-algebra
$\cal{B}$ of $A^\ZZ$.
When $\mu$ is $\sigma$-invariant the \emph{topological support} of $\mu$ is a subshift denoted by $S(\mu )$.
We shall say that the topological support is trivial if it is countable.
If $\alpha=(A_1,\dots,A_n)$ and $\beta=(B_1,\dots,B_m)$ are two partitions of
$X$, we denote by $\alpha \vee \beta$ the partition
$\{A_i\cap B_j, i=1,\dots,n,\; j=1,\dots,m\}$.
Let $T :X\to X$ be a measurable
continuous map on a compact set $X$.
The metric \emph{entropy} $h_\mu (T)$ of  $T$ is an isomorphism
invariant between two $\mu$-preserving transformations.
Let $H_\mu (\alpha ) =\sum_{A\in\alpha}\mu (A)\log \mu (A)$, where $\alpha$ is a finite partition of
$X$.
The entropy of the finite partition $\alpha$ is defined as
$h_\mu (T,\alpha )= \lim_{n\to\infty}1/nH_\mu (\vee_{i=0}^{n-1}T^{-i}\alpha )$
and the entropy of
$(X,T,\mu )$ as $h_\mu(T) = \sup_\alpha h_\mu (T,\alpha )$.

 A \emph{cellular automaton} is a continuous self-map $F$ on
$A^\ZZ$ commuting with the shift. The Curtis-Hedlund-Lyndon theorem \cite{He69}
states  that for every cellular automaton $F$ there exist an integer $r$ and a
{\em block map} $f:A^{2r+1}\mapsto A$ such that $F(x)_i=f(x_{i-r},\ldots ,x_i
,\ldots ,x_{i+r}).$
 The integer $r$ is called the \emph{radius} of the cellular automaton. If $X$
is a subshift of  $A^\ZZ$ and  $F(X)\subset X$, then the restriction of $F$
to $X$ determines a dynamical system $(X,F)$ called a cellular automaton  on
$X$.

\subsection{Equicontinuity, sensitivity  and expansiveness}

Let $F$ be a cellular automaton on  $A^\ZZ$.

\medskip
\hskip -.82 true cm{\bf Definition 1 (Equicontinuity)}
A point $x\in A^\ZZ$ is called an equicontinuous point (or  Lyapunov stable)
if for all  $\epsilon >0$, there exists  $\eta >0$ such that
$$d(x,y)\le \eta \  \Longrightarrow \  \forall i>0, \, d(T^i(x),T^i(y))\le
\epsilon.$$

\hskip -.82 true cm{\bf Definition 2 (Sensitivity)}
The automaton $(A^\ZZ ,F)$ is sensitive to initial conditions (or sensitive)  if
there exists a real number $\epsilon >0$ such that
$$
\forall x \in A^\ZZ, \, \forall \delta >0 ,\,
\exists y \in A^\ZZ,\, \mbox{  } d(x,y) \leq \delta, \mbox{ } \exists n\in\NN ,\,\, d(F^n(x), F^n(y))
\geq \epsilon.
$$

\medskip

The next definition appears in \cite{Gil87} for a Bernouilli measure.
\medskip

\hskip -.82 true cm{\bf Definition 3 ($\mu$-Expansiveness)}
The automaton $(A^\ZZ ,F)$ is $\mu$-expansive if there exists a real number $\epsilon >0$ such that
for all $x$ in $A^\ZZ$ one has
$$
\mu \left(\{y\in X:\, \forall i\in\NN,\, d(F^i(x),F^i(y))\le\epsilon \}\right)=0.
$$
\medskip

Notice that in \cite{Gil87} Gilman gives a classification of cellular automata based on
 the  $\mu$-expansiveness  and the $\mu$-equicontinuity classes.


\medskip

\hskip -.82 true cm{\bf Definition 4 (Positive Expansiveness)}
  The automaton $(A^\ZZ ,F)$ is positively exp\-an\-sive if there exists a real number $\epsilon >0$ such that,
\[
\forall (x,y)\in (A^\ZZ)^2, \,  x\neq y \, ,
\exists n\in\NN \mbox{ such that }d(F^n(x),F^n(y))\ge \epsilon .
\]
\medskip

K\r{u}rka \cite{Ku94} shows that, for cellular automata,
sensitivity is equivalent to the absence of  equicontinuous points.
\medskip

\subsection{Lyapunov exponents}
For all $x \in A^\ZZ$, the sets
$$
        W_s^+(x)=\{y\in A^\ZZ:\,\forall i\ge s, \, y_i=x_i \},\quad
        W_s^-(x)=\{y\in A^\ZZ:\,\forall i\le s, \, y_i=x_i\},
$$
are called \emph{right and left set of all the perturbations of $x$}, respectively.

\paragraph{}
For all integer $n$, consider the
smallest ``distance" in terms of configurations
at which a perturbation will not be able to influence the
$n$ first iterations of the automaton:
\begin{eqnarray}
        I^-_n(x)&=&\min\{ s\in\NN:\, \forall 1\le i\le n,\,
                F^i(W^-_{s}(x))\subset W_0^-(F^i(x))\}, \label{Im}\\
        I^+_n(x)&=&\min\{ s\in\NN:\, \forall 1\le i\le n,\,
                F^i(W^+_{-s}(x))\subset W_0^+(F^i(x))\}.
                \nonumber
\end{eqnarray}
We can now define the pointwise Lyapunov exponents  by
$$
\lambda^+(x)=\lim_{n\to\infty}\inf \frac{I^+_n(x)}{n}, \quad
\lambda^-(x)=\lim_{n\to\infty}\inf \frac{I^-_n(x)}{n}.
$$
For a given configuration $x$, $\lambda^+(x)$ and $\lambda^-(x)$
represent the speed to which the left and right faster perturbations
propagate.

\medskip

\hskip -.82 true cm{\bf Definition 5 (Lyapunov Exponents)}
For a $\mu$ shift-invariant measure on $A^\ZZ$,
we call average Lyapunov exponents of the automaton $(A^\ZZ ,F,\mu)$,
the constants
\begin{equation}\label{def.Lyap}
I^+_\mu=\liminf_{n\to\infty}\frac{I^+_{n,\mu}}{n},
\quad
I^-_\mu=\liminf_{n\to\infty}\frac{I^-_{n,\mu}}{n},
\end{equation}
where
\[
        I^+_{n,\mu}=\int_X I^+_n(x)d\mu (x),\quad
        I^-_{n,\mu}=\int_X I^-_n(x)d\mu(x).
\]

\begin{remark}
The sensitivity of the automaton $(A^\ZZ ,F,\mu)$ implies that for all $x\in A^\ZZ$,
$(I^+_n(x)+I^-_n(x))_{n\in\NN}$ goes to infinity.
\end{remark}

\section{Lyapunov exponents of positively expansive cellular automata}
\label{sectionpositivelyexpansive}

Similar versions of the next lemma appear in  \cite{Bl97} and \cite{Ku94}. The proof of similar
results in \cite{FMM98} using limsup is based on completely different arguments.
\begin{lem}\label{lem1}
Let $F$ be a positively expansive
CA with radius $r$ acting on a $F$-invariant subshift $X\subset A^\ZZ$.
There exists a positive integer $N^+$ such that for all
$x$ and $y$ in  $X$ that verify  $x(-\infty,-r-1)=y(-\infty, -r-1)$ and
 $F^m(x)(-r,r)=F^m(y)(-r,r)$
for all integers
 $0< m\le N^+$, we have  $x(r,2r)=y(r,2r)$.
\end{lem}
\paragraph{Proof}
 Let  $B_n$ be the subset of $(x,y)\in X \times X$ such
 that $x(-\infty,-r-1)=y(-\infty, -r-1)$, $x(r,2r)\neq y(r,2r)$
and $F^m(x)(-r,r)=F^m(y)(-r,r)$ for all $m<n$.
Each $B_n$ is closed and $B_{n+1}\subset B_n$.
Positive expansiveness of $F$ implies $\lim_{n\to\infty}B_n =
\emptyset$ (see \cite{Bl97}). Since $X$ is a compact set, there
is  a positive integer $N^+$ such that  $B_{N^+}=\emptyset$.
\hfill$\Box$

\begin{pro} \label{positiveexpansive}
For a positively expansive
CA acting on a bilateral subshift $X$,
there is a constant $\Lambda >0$ such that for all $x\in X$,
$\lambda^{\pm}(x)\ge \Lambda$.
\end{pro}

\paragraph{Proof}
We give the proof for $\lambda^-(x)$ only, the proof for $\lambda^{+}(x)$ being
similar.

Let $r$ be the radius of the automaton.
According to Lemma \ref{lem1},
 for any point $x\in X$ we obtain
 that  if $y\in W^-_{-1}(x)$ is such that
 $F^i(y)(-r,r)=F^i(x)(-r,r)$ ($\forall\; 1\le i\le N^+$) then $y$ must be in $W_{r}^-(x)\subset W^-_0(x)$.
 From the definition of $I_{N^+}^-(x)$ in (\ref{Im}), this implies that $I_{N^+}^-(x)\ge 2r$.
Lemma \ref{lem1} applied $N^+$  times  implies
that for each $0\le i\le N^+$,
$
F^n\left(F^i(x) \right)(-r,r)=F^n\left(F^i(y) \right)(-r,r)
$
for all $0\le n\le N^+$.
It follows that
$F^i(x)(r,2r)=F^i(y)(r,2r)$ ($\forall \; 0\le i\le N^+$).
Using Lemma \ref{lem1} once more and shifting $r$ coordinates of $x$ and $y$
yields $\sigma^r(x)(r,2r)=\sigma^r (y)(r,2r)\Rightarrow x(2r,3r)=y(2r,3r)$ so that
$I_{2N^+}^-(x)\ge 3r$.
Hence, for each integer $t\ge 1$, using Lemma \ref{lem1},
$N^+(t-1)!+1$ times  yields
$x\left(tr,(t+1)r\right)=y\left(tr,(t+1)r\right)$
and therefore  $I_{tN^+}(x)\ge (t+1)r$.
Hence for all $n\ge N^+$ and all $x\in X$,
$I_{n}^-(x)\ge (\frac{n}{N^+}+1)r$,
so that
$$
\lambda^-(x)=\lim_{n\to\infty}\inf \frac{I^-_n(x)}{n}\ge \frac{r}{N^+}.
$$
\hfill$\Box$


\section{Lyapunov Exponents and Entropy}

Let $F$ be a cellular automaton acting on a shift space $A^\ZZ$ and
let $\mu$ be a $\sigma$-ergodic and $F$-invariant probability measure.
According to the inequality
$$
        h_\mu (F)\le h_\mu (\sigma )(I^+_\mu+I^-_\mu)
$$
proved in \cite[Theorem 5.1]{ti2000}, one has
$I^+_\mu +I^-_\mu =0\Rightarrow h_\mu(F)=0$.
Here we extend this result to the case of a $\sigma$ and $F$-invariant measure on
a sensitive cellular automaton.

\begin{pro}\label{i=0}
If $F$ is a sensitive cellular automaton and $\mu$
a shift and $F$-invariant measure,
$I^+_\mu +I^-_\mu =0\,\Rightarrow\, h_\mu(F)=0.$
\end{pro}

\paragraph{Proof}
Let $\alpha$ be a finite partition of $A^\ZZ$ and
$\alpha_{n}^m (x)$ be the element of the partition
$\alpha\vee\sigma^{-1}\alpha\vee \dots \sigma^{-n+1}\alpha\vee \sigma^1\alpha
\dots\vee\sigma^m\alpha$ which
contains $x$.
Using \cite[Eq.~(8)]{ti2000} we see that, for all finite partitions $\alpha$,
\begin{equation}\label{ineq.hmu}
h_\mu (F,\alpha )\le \int_{A^\ZZ}\liminf_{n\to\infty}
      \frac{-\log\mu (\alpha^{I_n^+(x)}_{I_n^-(x)}(x)) }{I_n^+(x)+I_n^-(x)}
\times\frac{I_n^+(x)+I_n^-(x)}{n} d\mu(x).
\end{equation}
Suppose that
$
I^+_\mu +I^-_\mu=\liminf_{n\to\infty}\int_X n^{-1}\bigl(I_n^+(x)+I_n^-(x)\bigr)d\mu(x) =0$.
From Fatou's lemma we have $\int_X\liminf_{n\to\infty}
n^{-1}(I_n^+(x)+I_n^-(x))d\mu(x)=0$.
Since $n^{-1}(I_n^+(x)+I_n^-(x))$ is always a positive or null rational, there exists a set
$S\subset A^\ZZ$ of full measure such that $\forall x\in S$ we have
$\liminf_{n\to\infty}n^{-1}(I_n^+(x)+I_n^-(x))=0$.
Since $F$ is sensitive, for all points $x\in A^\ZZ$, we have $\lim_{n\to\infty}I_n^+(x)+I_n^-(x)=+\infty$
(see \cite{ti2000})
 and  the Shannon-McMillan-Breiman theorem (in the extended case of
 $\ZZ$ actions see \cite{Ki})
tells us that
$$
        \int_{A^\ZZ} \liminf_{n\to\infty}-
        \frac{\log\mu (\alpha_{I_n^-(x)}^{I_n^+(x)}(x)) }{I_n^+(x)+I_n^-(x)}d\mu =
h_\mu (\sigma,\alpha).
$$
Since for all $n$ and $x$, $-\log\mu
(\alpha_{I_n^-(x)}^{I_n^+(x)}(x))>0$, we deduce that  for all $\epsilon >0$
there is an  integer $M_\epsilon>0$ and a set $S_\epsilon\subset S$
with $\mu (S_\epsilon )>1-\epsilon$ such that
for all $x\in S_\epsilon$,
$$0\le\liminf_{n\to\infty}\frac{-\log\mu (\alpha_{I_n^-(x)}^{I_n^+(x)}(x))}{I_n^+(x)+I_n^-(x)}
\le M_\epsilon.
$$
For all $x\in S_\epsilon$ we obtain
$$
\phi(x) := \liminf_{n\to\infty}\frac{-\log\mu (\alpha_{I_n^-(x)}^{I_n^+(x)}(x)) }
{I_n^+(x)+I_n^-(x)}\times \frac{I_n^+(x)+I_n^-(x)}{n}=0,
$$
which implies  $\int_{S_\epsilon}\phi(x)d\mu(x) =0$.
Using the monotone convergence theorem
we deduce $\int_{A^\ZZ} \phi(x)d\mu(x) =0$.
It then follows from (\ref{ineq.hmu}) that
$h_\mu (F)=\sup_{\alpha} h_\mu(F,\alpha)=0$.
\hfill$\Box$

\section{The cellular automaton and its natural factor}

We define a cellular automaton for which the dynamic on a
 set of full measure is similar to  the "counters dynamic"  described in the introduction.
The unbounded size counters are "simulated" by the finite configurations in the
interval between two special letters "{\it E}". We will refer to these
special symbols as ``emitters''.
Between two $E$'s,
the dynamic of a counter is replaced by an odometer with overflow transmission to the right.
We add "$2$" and "$3$" to $\{0;1\}$ in the set of digits in order to have the sensitive dynamic
of counters with overflows transmission.
The states $2$ and $3$ are interpreted as ``$0$ + an overflow to be sent'',
and ``$1$ + an overflow to be sent''. Using this trick, we let
overflow move only one site per time unit. Notice that $3$ is necessary because it may happen that a counter is
increased by $2$ units in one time unit.

\subsection{The cellular automaton}

We define  a  cellular automaton $F$  from  $A^\ZZ$ to $A^\ZZ$ with $A=\{0;1;2;3;E\}$.
This automaton is the composition  $F = F_d \circ F_p$ of two cellular automata $F_d$
and $F_p$. The main automaton $F_d$ is defined by the local rule  $f_d$
\begin{eqnarray}\label{def-fd}
f_d(x_{i-2}x_{i-1}x_i)
&=& \mbox{  } {\bf 1}_{E}(x_i) x_i
  +  {\bf 1}_{\OE}(x_i)\left( x_i - 2\times{\bf 1}_{\{2, 3\}}(x_i)+{\bf 1}_{\{2,3\}}(x_{i-1})\right) \nonumber \\
&& {}\qquad+  {\bf 1}_{\OE}(x_i){\bf 1}_{E}(x_{i-1}) \left(1+{\bf 1}_{\{2\}}(x_{i-2})
\right),
\end{eqnarray}
where $\OE=A\setminus\{E\}$ and,
$$
   {\bf 1}_{S}(x_i)=\cases{1 & if $x_i\in S$,\cr
                        0 & otherwise.\cr
                        }
$$
The automaton $F_p$ is a ``projection'' on the subshift of finite type
made of sequences having at least three digits between two "E" which
is left invariant by $F_d$. Its role is simply to
restrict the dynamics to this subshift. It can be defined
by the local rule $f_p$
\begin{equation}\label{def-fp}
        f_p(x_{i-3},\dots,x_i,\dots x_{i+3})=
        {\bf 1}_{\OE}(x_i) x_i+
        {\bf 1}_{E}(x_i) x_i\times\prod_{j=-3\atop j\ne 0}^3{\bf
          1}_{\OE}(x_{i+j}).
\end{equation}
The dynamic of $F$ is illustrated in Figure \ref{fig1} for a particular
configuration.
The projection $f_p$ prevents the dynamic of $F$ to have equicontinuous points (points with
"blocking words"$EE$") and simplifies the relationship between the cellular automaton and the model.
The non-surjective cellular automaton acts surjectively on its $\omega$-limit space $X$ which is a
 non finite type subshift with a minimal distance of 3 digits between two "E". By definition,
 we have $X=\lim_{n\to\infty}\cap_{i=1}^n F^i(A^\ZZ)$.
The set $X$ is rather complicated. We do not want to give a
complete description. Some basic remarks may be useful for a
better understanding of the results. Note that the $E$s do not
change after the first iteration and that $x_i = 3$ implies
$x_{i-1} = E$. We can show that the word $222$ does not appear after
the second iteration and that $F^i(x)(k,k+1) = 22$ only if
$F^{i-1} (x) (k-2\colon k+1)\in\{2E21,0E31,1E31,2E31\}$. According
to the definition (\ref{def-fd}),  the evolution of finite
configurations without emitter "E" leads to sequences which
contains only the digits  "1" and "0".  This is the dynamic of the
emitter "E" with the overflows  crossing the "$E$s" that maintain
and move the letters "2" and "3". There is at most two letters "2"
between two consecutive letters "E".
 A typical word between two "E" is of the form $E3uE$, $E2uE$,
 $E3u_12u_2E$ or  $E2u_12u_2E$, where $u$, $u_1$ and
 $u_2$ are finite sequences of  letters "0" and "1" (the words $u_1$
 and $u_2$ can be empty).
Notice that  all  possible words  $u$, $u_1$, $u_2$ do not appear in $X$. For example, the
word $E200E$ does not belong to the language of $X$.

As we want to study the dynamic on finite (but unbounded)
counters, we define the set $\Omega \subset X$ of configurations
with infinitely many $E$ in both directions. This non compact set
is obviously invariant by the dynamics. We are going to define a
semi conjugacy between $(\Omega,F)$ and the model in the next
section.

\begin{figure}
\begin{center}
\begin{equation}\label{ex1}
\begin{array}{rcccccccccccc}
x =         \ldots & 0 & E & 1 & 1 & 0 & E & 0 & 2 & 2 & 2 & E& \ldots \\
F(x) =      \ldots & 0 & E & 2 & 1 & 0 & E & 1 & 0 & 1 & 1 & E& \ldots\\
F^{2}(x) =  \ldots & 0 & E & 1 & 2 & 0 & E & 2 & 0 & 1 & 1 & E& \ldots\\
F^{3}(x) =  \ldots & 0 & E & 2 & 0 & 1 & E & 1 & 1 & 1 & 1 & E& \ldots\\
F^{4}(x) =  \ldots & 0 & E & 1 & 1 & 1 & E & 2 & 1 & 1 & 1 & E& \ldots\\
F^{5}(x) =  \ldots & 0 & E & 2 & 1 & 1 & E & 1 & 2 & 1 & 1 & E& \ldots\\
F^{6}(x) =  \ldots & 0 & E & 1 & 2 & 1 & E & 2 & 0 & 2 & 1 & E& \ldots\\
F^{7}(x) =  \ldots & 0 & E & 2 & 0 & 2 & E & 1 & 1 & 0 & 2 & E& \ldots\\
F^{8}(x) =  \ldots & 0 & E & 1 & 1 & 0 & E & 3 & 1 & 0 & 0 & E& \ldots\\
F^{9}(x) =  \ldots & 0 & E & 2 & 1 & 0 & E & 2 & 2 & 0 & 0 & E& \ldots\\
\end{array}
\end{equation}
\caption{An illustration of the dynamic of $F$ defined by
         \ref{def-fp}--\ref{def-fd} on the configuration
        $x$, assuming that $x$ is preceded by enough $0$.}
\label{fig1}
\end{center}
\end{figure}

\subsection{The natural factor}

In order to make more intuitive the study of the dynamic of $F$ and to define (see Section 7) a natural
measure,
we introduce the  projection of this CA
which is a continuous  dynamical system that commutes with an infinite state 1-dimensional shift.


The word between two consecutive $E$s can be seen as a ``counter" that
overflows onto its right neighbour when it is full.
At each time step $E$  "emits" $1$ on its right
except when the counter on its left overflows: in this case there is a carry
of $1$  so the  $E$ "emits" $2$ on its right.
$$
\begin{array}{ccccccccccccccc}
x &=&         \ldots & 0 & \overbrace{E}^{\mbox{\tiny emitter $i$}} &
             \ 0 & 0 & 2 & \overbrace{E}^{\mbox{\tiny emitter $i+1$}}
             & \ 1 & 1 & 0 & 0 & \overbrace{E}^{\mbox{\tiny emitter $i+2$}} & \ldots \\
F(x) &=& \ldots & 0 & E &
            \multicolumn{3}{c}{\underbrace{\begin{array}{ccc}
            1 & 0 & 0\end{array}}_{\mbox{\tiny counter $i$}}}
            & E &
            \multicolumn{4}{c}{\underbrace{\begin{array}{cccc}
            3 & 1 & 0 & 0\end{array}}_{\mbox{\tiny counter $i+1$}}}
            & E& \ldots\\
\end{array}
$$
In what follows we call {\em counter} a triple
$(l,c,r)$, where
$l$ is the number of digits of the counter, $c$ is its state and
$r$ is the overflow position in the counter.
In Figure 1, in the first counter the countdown starts in $F^5(x)$ when the "2" is followed by "11".
This "2" can propagate at speed one to the next emitter "E".

Recall that $\Omega \subset X$ is a set
of configurations with
infinitely many $E$s in both directions and has at least three digits
between two $E$s.
Define the sequence
$(s_j)_{j \in \ZZ}$ of the positions of the $E$s in $x\in\Omega$ as follows:
\begin{eqnarray*}
        s_0(x) &=& \sup{\{ i \leq 0 \; : \; x_i = E\} },\\
        s_{j+1}(x) &=&  \inf{\{ i > s_j(x) \; : \; x_i = E\}}
        \quad \mbox{for $j \geq 0$,}\\
        s_{j}(x) &=& \sup{\{ i < s_{j+1}(x) \; : \; x_i = E\} } \quad
        \mbox{for $j<0$.}
\end{eqnarray*}
Denote by $u=(u_i)_{i \in \ZZ} =   (l_i, c_i,r_i)_{i \in \ZZ}$
a bi-infinite sequence of counters. Let $B=\NN^3$ and $\sigma_B$ be the shift on $B^\ZZ$.
We are going to define a function $\varphi$ from $\Omega\to B^\ZZ$.
We set for all $i \in \ZZ$, $l_i(x) = s_{i+1}(x) - s_i(x) - 1$ and define
$d_i(x) = \sum_{j=s_i +1}^{s_{i+1} -1} x_j 2^{j-s_i-1}$.
We denote by $\overline{c}_i = 2^{l_i}$ the period of the counter
$(l_i,c_i,r_i)$.

For each $x\in\Omega$ and $i\in\ZZ$ we set  $c_i(x) = d_i(x) \hbox{ modulo } \overline{c}_i(x)$
and we write
$$
r_i(x) =
\left\{
\begin{array}{ll}
l_i+1-\max{\{ j \in \{s_i +1, \ldots, s_{i+1} -1 \} \; : \; x_j >1 \}
} + s_i& \hbox{ if } d_i(x) \geq
 \overline{c}_i(x)\\0 & \hbox{ otherwise.}
\end{array}
\right.
$$
For each $x\in\Omega$ we can define $\varphi (x) =(l_i(x),c_i(x),r_i(x))_{i\in\ZZ}$.
Remark that since $r_i(x)\le l_i(x)$ and $c_i(x)\le 2^{l_i(x)}$,
the set $\varphi (\Omega)$ is a strict subset $(\NN^3)^\ZZ$ .

\begin{figure}
\begin{center}
$$\begin{array}{rccc}
u = \ldots  & (3,3,0) & (4,0,0) & \ldots \\
H(u) =  \ldots  & (3,4,0) & (4,1,0) & \ldots \\
H^{2}(u) = \ldots  & (3,5,0) & (4,2,0) & \ldots \\
H^{3}(u) = \ldots  & (3,6,0) & (4,3,0) & \ldots \\
H^{4}(u) =  \ldots  & (3,7,0) & (4,4,0) & \ldots \\
H^{5}(u) =  \ldots  & (3,0,3) & (4,5,0) & \ldots \\
H^{6}(u) =  \ldots  & (3,1,2) & (4,6,0) & \ldots \\
H^{7}(u) =  \ldots  & (3,2,1) & (4,7,0) & \ldots \\
H^{8}(u) =  \ldots  & (3,3,0) & (4,9,0) & \ldots \\
H^{9}(u) =  \ldots  & (3,4,0) & (4,10,0) & \ldots \\
\end{array}
$$
\caption{The dynamic in Figure \ref{fig1} for the natural factor.}
\label{fig2}
\end{center}
\end{figure}

\paragraph{}
On $\varphi (\Omega)$, we define a dynamic on the counters through a local function.
First we give a rule for incrementation of the counters.
For  $a=1$ or $2$, we set
\[
\left\{
\begin{array}{lllll}
(l_i,c_i, 0) + a & = & (l_i,c_i+a, 0) &\hbox{ if } c_i < \overline{c}_i -a & (R1)\\
(l_i,c_i, 0) + a &=& (l_i,c_i+a - \overline{c_i}\,,l_i)&\hbox{ if } c_i+a \ge \overline{c}_i&(R2)\\
(l_i,c_i, r_i) + a &= &(l_i,c_i + a, r_i -1) &\hbox{ if } r_i >0 &(R3).
\end{array}
\right.
\]
We define the local map $h$ on $\NN^3 \times\NN^3$  by
\[
h(u_{i-1}u_i)= u_i+ (1+ \carac{r_{i-1}=1}(u_{i-1})),
\]
where the addition must be understood following the incrementation
procedure above, with $ a = 1+ \carac{r_{i-1}=1}(u_{i-1})$.
Let $H$ be the global function on $(\NN^3)^\ZZ$.
This is a ``cellular automaton on a countable alphabet''.

\paragraph{}
Note that, the $(l_i)_{i \in \ZZ} $ do not move under iterations.
At each step, the counter $c_i$ is increased by  $a$ modulo
$\overline{c}_i$ ($a=1$ in general, while $a=2$ if counter $i-1$
``emits'' an overflow).
When $c_i$ has made a complete turn $r_i$ starts to count  down $l_i$,
$l_{i}-1\ldots 1, 0 $ ; after $l_i$ steps  $r_i$ reaches $0$,
indicating that (in the automaton) the overflow has reached its position.

\begin{remark}
Notice that if $2 l < 2^l$, we cannot have $c_i = 2^l-a$ and $r_i>0$ since
$r_i$ is back to $0$ before $c_{i-1}$  completes a new turn. This technical detail is the reason
why we impose a minimal
distance $3$ (more than the distance one required for the sensitivity condition) between two successive $E$.
\end{remark}

\subsection{Semi conjugacy}

\begin{pro}
We have the following semi conjugacy,
$$
\begin{array}{lcl}
&F&\\
\Omega & \longrightarrow & \Omega\\
\downarrow \varphi & & \downarrow \varphi\\
&H&\\
\varphi (\Omega ) & \longrightarrow &\varphi (\Omega )\\
\end{array}
$$
with
$$
\varphi \circ F = H \circ \varphi.
$$
\end{pro}

\paragraph{Proof}

Let $x \in \Omega$. Denote $\varphi(x) = (l_i,c_i,r_i)_{i \in \ZZ}$,
$x'= F(x)$ and $\varphi(x') = (l'_i,c'_i,r'_i)_{i \in \ZZ}$. We have to prove that
 $(l_i, c_i, r_i) + 1 + {\bf 1}_{\{r_{i-1}=1\}} =  (l'_i,c'_i,r'_i)$ where the addition
 satisfies the rules $R_1, R_2, R_3$.

First, we recall that the $Es$ do not move so that $l'_i = l_i$.

Consider the first digit after the ith emitter
$E$ :$x'_{s_i+1}=x_{s_i+1}-2\times {\bf 1}_{\{2,3\}}(x_{s_{i}+1})+1+{\bf 1}_{\{2\}}(x_{s_i-1})$.
Clearly we have $r_{i-1}=1$ if and only if $x_{s_i-1}=2$
 so $x'_{s_i+1}=x_{s_i+1}+1+{\bf 1}_{\{r_{i-1}=1\}}-2\times {\bf 1}_{\{2,3\}}(x_{s_{i}+1})$.
For all $s_i+2\le j\le s_{i+1}-1$,  $x'_j=x_j-2\times {\bf 1}_{\{2\}}(x_j)+{\bf 1}_{\{2,3\}}(x_{j-1})$.
Since $d_i(x) = \sum_{j=s_i +1}^{s_{i+1} -1} x_j 2^{j-s_i-1}$, if $x_{s_{i+1}-1}\neq 2$ then
 $d_i(x')=d_i(x)+1+{\bf 1}_{\{r_{i-1}=1\}}$ and if $x_{s_{i+1}-1}=2$ then
 $d_i(x')=d_i(x)+1+{\bf 1}_{\{r_{i-1}=1\}}-2\times 2^{-l_i}=d_i(x)+1+{\bf 1}_{\{r_{i-1}=1\}}-\overline{c_i}$.
 As $c'_i=d_i(x')\, \mbox{mod} \;\overline{c_i}\;$ then for all $x\in\Omega$
  one has  $c'_i=c_i+1+{\bf 1}_{\{r_{i-1}=1\}}\, \mbox{mod} \;\overline{c_i}$.

 It remains to understand the evolution of the overflow $r_i$.
 First, notice that if $d_i(x')=\overline{c_i}=2^{l_i}$ then $x'(s_{i},s_{i+1})=E21^{(l_{i}-1)}E$
 and if $d_i(x')=\overline{c_i}+1=2^{l_i}+1$ then $x'(s_{i},s_{i+1})=E31^{(l_{i}-1)}E$.
 After $l_i$ iteration of $F$, the configurations $E21^{(l_{i}-1)}E$ and $E31^{(l_{i}-1)}E$ have
 the form $Ew2E$.
 The maximum value taken by $d_i(x)$ is when $x(s_{i},s_{i+1})=Eu2E$ where $d_i(z)<2l_i$
 if $z(s_{i},s_{i+1})=Eu0E$. As noted in Remark 2, the counters,
which have at least a size of 3, have not the time to make a complete turn during the
 countdown ($2l_i<\overline{c_i}$) which implies that $d_i(x)\le \overline{c_i}+2l_i <2\overline{c_i}$.
Since each counters can not receive
 an overflow at each iteration then $d_i(x)<2\overline{c_i}-2$
 (when  $l_i\ge 3$, $l_{i-1}\ge 3$, $\theta\times l_i< \overline{c_i}-2=2^{l_i}-2$ where $\theta <2$ ).
 Clearly $r_i=0$ if and only if $d_i=c_i<\overline{c_i}$  (addition rule $R_1$).
As $d_i<2\overline{c_i}-2$, if  $c_i=\overline{c_i}-2=d_i(x)$ and $x_{s_i-1}=2$ or $c_i=\overline{c_i}-1=d_i(x)$
($r_i=0$)
then $x'(s_{i},s_{i+1})=E21^{(l_{i}-1)}E$ or $x'(s_{i},s_{i+1})=E31^{(l_{i}-1)}E$ which implies
that $r'_i=l_i$  (addition rule $R_2$).

  Now remark that if $r_i>0$ then
$x(s_{i},s_{i+1})=Eu21^kE$ with $0\le k\le l_i-1$ and $u$ is a finite sequence of letters
"0", "1", "2" or "3". Using the local rule of $F$, we obtain that the letter "2" move to the right of one
coordinate which implies that $r'_i=r_i-1$ (addition rule $R_3$).

\hfill$\Box$
\begin{remark}
\label{invertible}
We remark that $\varphi$ is not injective.
 Consider the subset of $X$ defined by
$$\Omega^* = \Omega \cap \{ x \in \{0,1,E\}^\ZZ \; : \; x_0 = E \}. $$
It is clear that $\varphi$ is one to one between $\Omega^*$ and $\varphi (\Omega)$ since the origin is fixed and
there is only one way to write the counters with $0$ and $1$. We will use this set,
keeping in mind that it is not invariant for the cellular automaton $F$.
\end{remark}

\begin{remark}
If $x \in A^\ZZ$, 
we use
$c_i(x)$ to denote $c_i(\varphi(x))$ and $l_i(x)$ instead of $l_i(\varphi (x))$.
 This should not yield any confusion.
To take the dynamics into account, we write $c_i^t(x) =
c_i(\varphi(F^t(x)))= c_i(H^t(\varphi(x)))$ and similarly
 for $l_i(x)$ and $r_i(x)$. Note that $l^t_i(x) =
l_i(x)$ for all $t$.
\end{remark}

\subsection{Limit  periods}

The natural period of the counter $i$, is its number of states, say $\overline{c}_i$.
We introduce the notion of {\em real period} or
{\em asymptotic period} of a counter which is, roughly speaking,
the time mean of the successive observed periods.
It can be formally defined as the inverse of the
number of overflow emitted by the counter (to the right)
per unit time. More precisely, we define $N_i(x)$, the inverse of the real period $p_i(x)$ for the counter number $i$ in
$x$ by
$$N_i(x) = \lim_{t \to +\infty} \frac{1}{t} \sum_{k=0}^{t}
\carac{1}(r_i^t(x)).$$
\begin{lem}\label{limperiod}
\label{periode}
The real period $p_i(x)$ exists and,
$N_i(x) =\frac{1}{p_i(x)}= \sum_{k=i}^{-\infty}   2^{ - \sum_{j=i}^k l_j(x)}. $
\end{lem}

\paragraph{Proof}
Let $n_i^t$ be the number of overflows emitted by the counter $i$ in $x$ before time $t$.
$$n_i^t(x) =  \sum_{k=0}^{t}  \carac{1}(r_i^t(x)).$$
The number of turns per unit time is the limit when $t \to + \infty$ of $n_{i}^t/t$ if it
exists. We can always consider the limsup $N_i^+$ and the  liminf $N_i^-$
of these sequences.

After $t$ iterations, the  counter indexed by $i$  has been incremented  by  $t$
(one at each time step) plus the number of overflow ``received'' from the counter $(i-1)$.
We are looking for a recursive relationship between $n^t_i$ and
$n^t_{i-1}$.
Some information is missing about the delays in the ``overflow
transmission'', but
we can give upper and lower bounds.

The number of overflows emitted by counter $i$ at time $t$ is essentially given by
its initial position + its ``effective increase'' - the number of overflows delayed,
divided by the size $\overline{c}_i$ of the counter. The delay in the overflow transmission
is at least $l_i$. The initial state is at most $\overline{c}_i$.
Hence, we have,
\[
n_i^t(x) <\frac{\overline{c}_i(x)+t+n_{i-1}^t(x)}{\overline{c}_i(x)}
\hskip .5 true cm \mbox { and, } \hskip .5 true cm
n_{i}^t(x)>\frac{t+n_{i-1}^t(x) -l_i(x)}{\overline{c}_i(x)}.
\]

So for the limsup $N_i^+$ and liminf $N_i^-$, we obtain
$ N^\pm_i=\frac{1}{\overline{c}_i}(1+N^\pm_{i-1})$.
Remark that for all $t\in\NN$ we have
$N_i^+-N_i^-=\frac{1}{\overline{c}_i}(N_{i-1}^+-N_{i-1}^-)=(\frac{1}{\overline{c}_i})(N_{i-t}^+-N_{i-t}^-)$.
Since $\frac{1}{\overline{c_i}}\le N_i^-\le N_i^+\le \frac{2}{\overline{c_i}}$,
the limit called $N_i$ exists and finally we have
\begin{equation}
N_i(x) =  \sum_{k=i}^{-\infty}  \prod_{j=i}^k    \overline{c}_j^{-1}(x)
=  \sum_{k=i}^{-\infty}   2^{ - \sum_{j=i}^k l_j(x)}.
\end{equation}
\hfill$\Box$
\begin{remark}
The series above are lesser than the convergent geometric series $(\sum_{k=1}^n$ $2^{-3k})_{n\in\NN}$
 since $l_i$ is always greater than 3.
Note that in the constant case, $l_i = L$, the limit confirm the intuition because  the period is
$$p = \frac{1}{N_i} = \frac{1}{\sum_{k=1}^{+\infty} 2^{-kL}}= 2^L - 1. $$
\end{remark}

\section{Sensitivity}

For the special cellular automata $F$, we say that a measure $\mu$ satisfies conditions
 (*) if for all $l\in (\NN)^\ZZ$ one has $\mu (A^\ZZ \setminus \Omega )=0$ and
 $\mu (\{x\in\Omega : \left(l_i(x)\right)_{i\in\ZZ}=l\})=0$.
 These ``natural'' conditions are satisfied by the invariant
 measure $\mu^F$ we consider (see Section 7, Remark 6).

\begin{pro}\label{sensitivity}
The automaton $F$ is sensitive to initial conditions. Moreover, it is  $\mu$-expansive
if $\mu$ satisfies condition (*).
\end{pro}

\paragraph{Proof}

Fix $\epsilon=2^{-2}=\frac{1}{4}$ as the sensitive and
$\mu$-expansive constant. When  $x\in\Omega$,  it is possible to
define $l(x)$ which is the sequence of the size of the counters
for $x$ and to use the model $(H, \varphi (\Omega ))$ to
understand the dynamic.

We can use Lemma \ref{periode}  to prove that for the model,  if
we modify  the negative coordinates of the sequence $l_i$ we also
modify the asymptotic behaviour of the counter at $0$. Such a
change in the asymptotic behaviour implies that at one moment, the
configuration at $0$ must be different. From the cellular
automaton side, we will show that a change of the real period of
the "central counter" will affect the sequences
$(F^t(x)(-1,1))_{t\in\NN}$ which is enough to prove the
sensitivity and $\mu$-expansiveness conditions.

For each $x\in \Omega$ with $x_0\neq E$,  consider the sequence $l_{[0,-\infty]}(x)=
l_{0}(x)l_{-1}(x)$ $\ldots  l_{-k}(x) \ldots$.
We claim that if  $l_{[0,-\infty]}(x)\neq l_{[0,-\infty]}(y)$ then  $N_0(x)\neq N_0(y)$.
Let $j$ be the first negative or null integer such that $l_j(x)\neq l_j(y)$.
By Lemma \ref{periode}, there exists a positive real $K_1=\sum_{k=0}^{j-1}2^{-\sum_{i=0}^{k}l_i(x)}$ such that
$$
N_0(x)=K_1+K_1 2^{-l_j(x)}+K_1 2^{-l_j(x)}(\sum_{k=j+1}^\infty 2^{-\sum_{i=j+1}^k l_i(x)})
$$
and
$$
N_0(y)=K_1+K_1 2^{-l_j(y)}+K_1 2^{-l_j(y)}(\sum_{k=j+1}^\infty 2^{-\sum_{i=j+1}^k l_i(y)}).
$$
So writing $K_2=K_1 2^{-l_j(x)}$ we obtain
$$
N_0(x)-N_0(y)>\frac{K_2}{2}-\frac{K_2}{2}(\sum_{k=j+1}^\infty 2^{-\sum_{i=j+1}^k l_i(y)}).
$$
As $\sum_{k=i+1}^\infty 2^{-\sum_{j=i+1}^k l_j(y)}$ is less than the geometric series
$\sum_{k=0}^\infty 2^{-3k}=\frac{1}{7}$ we prove  the claim.

If $x_0=E$, using the shift commutativity of $F$ we obtain that if
$l_{[-1,-\infty]}(x)\neq l_{[-1,-\infty]}(y)$ then  $N_{-1}(x)\neq N_{-1}(y)$.

\medskip

Remark that for each $x\in\Omega$ and $\delta =2^{-n}>0$, there exists $y\in\Omega$ such that
$y_i=x_i$ for $i\ge -n$ and $l(x)\neq l(y)$ (sensitivity conditions).
We are going to show that if
$l_i(x)\neq l_i(y)$ ($i<0$; $x\in\Omega$) then there exist some $t\in\NN$ such that
$F^t(x)(-1,1)\neq F^t(y)(-1,1)$.

First, we consider the case where $x_0 = E$. In this case, since
$N_{-1}(x) \neq  N_{-1}(y)$, there is $t$ such that $n^t_{-1}(x) \neq
n^t_{-1}(y)$. At least for the first such $t$, $r^t_{-1}(x) \neq
r^t_{-1}(y)$ since $r^t_{-1}(x)=1$ and $r^t_{-1}(y)\neq 1$.
But this exactly
means that $x^t_{-1} =2$ and $y^t_{-1} \neq 2$.
Hence $x^t_{-1} \neq y^t_{-1}$.

Next assume that $x_0 \neq  E$.
We have $N_0(x) > N_0(y)$. This implies that $n_0^{t}(x) - n_0^{t}(y)$
goes to infinity. Hence the difference $c^t_0(x) - c^t_0(y)$ (which can
move by $0$, $1$ or $-1$ at each step)  must take
(modulo $\overline{c}_0$)
all values between
$0$ and $\overline{c}_0-1$. In particular, at one time $t$,
the difference must be equal to  $2^{-T_0-1}$.
If $x_{-1} = E$, then $T_0 = -1$ which implies that $x^t_0 \neq
y^t_0$. Otherwise, in view of the conjugacy, it means that
$$\sum_{i=T_0 +1}^{T_{1} -1} x^t_i 2^{i-T_0-1} - \sum_{i=T_0 +1}^{T_{1} -1} y^t_i
2^{i-T_0-1} = 2^{-T_0-1} \mbox{(modulo $\overline{c}_0$)},$$ that is, either $x^t_0 \neq
y^t_0$ or $x^t_0 = y^t_0$ and $x^t_{-1} \neq y^t_{-1}$.
We have proved the sensitivity of $F$ for $x\in \Omega$.

\medskip
To show that $x\in\Omega$ satisfies the $\mu$-expansiveness condition,
 we need to prove that the set of points which have the same asymptotic period for the central counter
is a set a measure zero. For each point $x$ the set of all points $y$ such that
$F^t(x)(-1,1)=F^t(y)(-1,1)$ or $d(F^t(x),F^t(y))\le \frac{1}{4}$ ($t\in\NN$),
is denoted by $D(x,\frac{1}{4})$. Following the arguments of the proof of the sensitivity condition above,
we see that every change in the sequence of letters ``Es'' in the left coordinates will affect the
central coordinates after a while, so  we obtain that  $D(x,\frac{1}{4})\subset \{y : l_i(y)=l_i(x) : i< 0 \}$.
Since $\mu$ satisfies condition (*), we have  $\mu (D(x,\frac{1}{4}))=0$ which is the condition
required for the  $\mu$-expansiveness.

Now suppose that  $x\in \{0,1,2,3,E\}^\ZZ \setminus \Omega$. First
notice that if there is at
least one letter $E$ in the
left coordinates, the sequence $F^t(x)(-1,1)$ does not depend on the position or even the existence
of a letter $E$ in the right coordinates.
Recall that after one iteration, the word $22$ appears only directly after a letter $E$.
So when there is at least one letter $E$ in the negative coordinates of $x$, the arguments and the proof given
for $x\in\Omega$ still work.

If there is no letter $E$ in the negative coordinates of $x$, for any $\delta=2^{-n}$, we can consider
any $y\in\Omega$ such that $y_i=x_i$ if $-n\le i\le n$. For $x\in \{0,1,2\}^\ZZ$, the dynamic is only
given by the letter $2$ which move on sequences of $1$ (see Figure 1).  The sequence $F^t(x)(-1,1)$ can not behave like a
counter because it is an ultimately stationary sequence (after a letter ``2'' pass over a ``1'' it remains a ``0'' that can not
be changed)  and the real period of the central counter of $y$ is obviously strictly greater than 1.
Then  in this case again,  the sequences
 $F^t(x)(-1,1)$ and $F^t(y)(-1,1)$ will be different  after a while which satisfies the sensitivity condition.
Since the only points $z$ such that $F^t(z)(-1,1)=F^t(x)(-1,1)$ belong to the set of null
measure $\{0,1,2,3\}^\ZZ$ (as $\mu$ satisfies conditions (*)),  we obtain the $\mu$-expansiveness condition.
\hfill$\Box$

\section{Invariant measures}

\subsection {Invariant measure for the model}

We construct an invariant measure for the dynamics of the counters.
Firstly, let $\mulength$ denote a measure on $\NN \setminus \{0,1,2\}$ and secondly, let
$\muemetteur = \mulength^{\otimes \ZZ}$ be the product measure. To fix
ideas, we take $\mulength$ to be the geometric law of paramater
$\nu = \frac{2}{3}$ on $\NN$ conditionned to be larger than  $3$, i.e.,
$$\mulength(k) =
\left\{
\begin{array}{ll}
\frac{\nu^k}{\sum_{j>2} \nu^j} & \hbox{if } k\geq 3\\
0 &\hbox{if } k < 3.
\end{array}
\right.
$$
Notice that the expectation of $l_0$ is finite : $\sum_{l_0 >2}\nu_*(l_0)\times l_0
 < +\infty$.
We denote by $m_L$ the uniform measure on the finite set $\{0, \ldots, L-1\}$.
Given a two-sided sequence $l  = (l_i)_{i \in \ZZ}$, we define a
measure on $\NN^\ZZ$ supported on $\prod_{i \in \ZZ} \{0, \ldots,
2^{l_i}-1\}$, defining
$\mucompteurfix_{l} = \otimes_{i \in \ZZ} \muunif_{2^{l_i}}$, so that for all
$k_{0}, \ldots , k_{m}$ integers,
$$ \mucompteurfix_{\bf l}(\{c_{i} = k_0, \ldots , c_{i+m} = k_m\} ) =
\left\{
\begin{array}{ll}
2^{-\sum_{j=i}^{i+m} l_j} & \hbox{ if, }   \forall i\leq j \leq i+m, \, c_j < 2^{l_j}\\
0  & \hbox{otherwise}.
\end{array}
\right.
$$
We want this property for the counters to be preserved by the
dynamics.
But, for the overflow, we do not know a priori how the measure will
behave.
Next we construct an initial measure and iterate the sliding block code on the counters $H$.
For all ${l} = (l_i)_{i \in
  \ZZ}$, we set  $\nu_{l} = \otimes_{i \in \ZZ} \delta_{l_i}$,
and,
$\muoverflow = \otimes_{i \in \ZZ} \delta_0$, where $\delta_k$ denotes
the Dirac mass at integer $k$. We consider the measures
$$\mu^H_{l} = \nu_{l} \otimes \mucompteurfix_{l} \otimes
\muoverflow, $$
 which  give
mass $1$ to the event $\{ \forall i \in \ZZ, r_i= 0 \}$.
Then, we define,
$$\tilde \mu^H= \int_{\NN^\ZZ} \muglob_{ l}\;  \muemetteur(d{l}).$$
The sequence $ \frac{1}{n} \sum_{k=0}^{n-1}  \tilde \mu^H   \circ H^{-k}$ has
convergent subsequences. We choose one of these subsequences, $(n_i)$,
and write
$$
\muglob = \lim_{i\to\infty}\frac{1}{n_i}\sum_{k=0}^{n_i-1} \tilde
\mu^H  \circ H^{-k}.
$$

\subsection{Invariant measure for $F$}

Now we construct a shift and $F$-invariant measure for the cellular automaton space.
\paragraph{}
Recall that $\Omega^* = \Omega \cap \{ x \in \{0,1,E\}^\ZZ \; : \; x_0 = E \}$
(see Remark \ref{invertible}). We write, for all integer $k$, $\Omega^*_k =
\sigma^k \Omega^*$.
Let us define, for all  measurable subsets $I$
of  $A^\ZZ$,
$$\muunf   = \sum_{k=0}^{l_0-1}  \muglob_{l}\left(\varphi(I\cap
  \Omega_k^*)\right). $$
This measure distributes the mass on the $l_0$ points with the same image (in
the model) corresponding to the $l_0$ possible
shifts of  origin.
The total mass of this measure is $l_0$.
Since the expectation of $l$  is
$\bar{l}=\sum_{i=3}^\infty \muemetteur_*(l)\times l < \infty$ is finite
(if $\nu=\frac{2}{3};\; \overline{l}=5$), we can  define a probability measure
$$
\tilde{\mu}^F= \frac{1}{\overline{l}}  \int_{\NN^\ZZ} \mudeuxf\;  \muemetteur(d{l}).
$$
The measure $\tilde \mu^F$ is shift-invariant (see further) but it is
supported on a non $F$-invariant set.
In order to have a $F$-invariant measure we take an adherence value of
the Cesaro mean. We choose a convergent subsequence $(n_{i_j})$ of the sequence
$(n_i)$ defining $\mu^H$ and write
$$\mu^F = \lim_{j \to \infty} \frac{1}{n_{i_j}}\sum_{k=0}^{n_{i_j}-1}\tilde{\mu}^F \circ F^{-k}.$$

\begin{remark}
Since for the measure $\mu^F$, the length between two ``$Es$'' follows a geometric law, the measure
$\mu^F$ satisfies the conditions (*) defined in Section 6, and from Proposition 5, the automaton
 $F$ is $\mu^F$-expansive.
\end{remark}

\begin{lem}
\label{proprietesmesures}
The measure $\muglob$ is $\sigma_B$ and $H$-invariant.
The measure $\muf$ is  $\sigma$ and $F$-invariant.
For all measurable subsets $U$ of $\varphi (\Omega )$  such that $l_0$ is constant
on $U$,
$$
\mu^F(\varphi^{-1}(U))=
\frac{l_0}{\overline{l}} \muglob(U).
$$
\end{lem}
\paragraph{Proof}
The shift invariance of $\muglob$ follows from the fact that $\muemetteur$
is a product measure and that the dynamic commutes with the shift so
the Cesaro mean does not arm.
$H$ and $F$ invariance of $\muglob$ and $\muf$ follow from the standard argument
on Cesaro means.

\paragraph{}
Shift invariance of $\mu^F$ comes from a classical argument of Kakutani towers
because $\mu^H$ is essentially the induced measure of $\mu^F$ for the shift on
the set $\Omega^*$. We give some details.
We choose a measurable set $I$ such that $l_0$ is
constant on $I$.
Choose $l = (l_i)_{i \in \ZZ}$. Noticing that $\varphi(\sigma^{-1}(I
)\cap \Omega^*_k) = \varphi(I
\cap \Omega^*_k)$ as soon as $k\geq 1$,
a decomposition of $I$ such as
$ I = \cup_{k=0}^{l_0-1} I \cap \Omega^*_k$ yields
\begin{eqnarray*}
\muunf(\sigma^{-1} I)  &=&
\muunf(I)- \muunf(I \cap \Omega^*)
+ \muunf( \sigma^{-1}(I) \cap \Omega^*).
\end{eqnarray*}
We remark that
$$\varphi(\sigma^{-1}(I)
\cap \Omega^*) =  \sigma_B^{-1} (\varphi(I\cap \Omega^*))$$
so that
$$\muunf( \sigma^{-1}(I) \cap \Omega^*) =  \muglob_{l}\left(\sigma_{B}^{-1}
  (\varphi(I  \cap \Omega^* )) \right) . $$
Now we integrate with respect to $\bf l$ and use the
$\sigma^{-1}_B$-invariance
of $\tilde \mu^H$ to conclude that
$$\tilde \mu^F(\sigma^{-1} I) = \tilde \mu^F(I).$$
For a measurable set $I$, we decompose $I= \cup_{L \in \NN}
(I\cap \{x\in \Omega : l_0(x) = L\} ) $.
Since $F$ commutes with $\sigma$, the Cesaro mean does not affect the
shift
invariance. Hence, $\mu^F$ is shift invariant.

\paragraph{}
Now, we make explicit the relationship  between $\muglob$ and $\muf$.
Let $V \subset \Omega$ be an event
measurable with respect to the $\sigma$-algebra generated by
$(l_i,c_i ; i \in \ZZ)$. There is
an event $U \subset \varphi (\Omega )$ such that $V= \varphi^{-1}(U)$.
Write $V = \cup_{L \in \NN} V_L$,  where $V_L = V \cap \{x\in\Omega : l_0(x) = L\}$.
Since $\varphi(V_L \cap
\Omega^*_k) = \varphi(V_L) =: U_L$ for all $0 \leq k <L$, we can write
$$\muunf(V_L) = \sum_{k=0}^{L-1} \mu^H_l (\varphi(V_L \cap
\Omega^*_k)) = L \mu^H_l(\varphi(V_L)).$$
Integrating with respect to $l$ gives
$$\tilde \mu^F (V )=  \frac{1}{\overline{l}}  \sum_{L\in \NN} L  \tilde \mu^H  (U_L) \muemetteur (L).$$
 If $l_0$ is constant on $U$ then $L=l_0$ and  we get
$$\tilde  \mu^F (\varphi^{-1}(U))  =  \frac{L}{\overline{l}}  \tilde \mu^H  (U)$$
Since $F^{-1} (\varphi^{-1}(U)) = \varphi^{-1} (H^{-1}(U))$ and   the Cesaro means are taken along the same subsequences, we are
done.
\hfill$\Box$

\begin{lem}
\label{uniforme}
Fix $M \in \ZZ, M<0$ and a sequence $L= (L_{M}, \ldots,
L_0)$. Consider the cylinder $V_L = \{x \in \Omega : (l_M(x), \ldots,
l_0(x)) = L \} $.
For all $M<i<0$, and all $0 \leq k < 2^{L_i}$,
$$ \mu^F (\{x:c_{i}(x)=k\} | V_L )= 2^{-L_i}.
$$
\end{lem}
\paragraph{Proof}
Let ${\cal F}$ be the $\sigma$-algebra generated by $\{ (c_j,r_j)\vert j<0\}$
and $\{ l_j, j \in \ZZ\}$.
We claim that, almost surely,
$$\muglob (\{c_0 = k\} \;|\;{\cal F} ) =
\mucompteurfix_{l}(\{c_0 = k\})   = 2^{-l_0}. $$
Let $u, v \in \varphi (\Omega )$.
If $u_i = v_i$ for all $i<0$, and $c_0(v) = c_0(u) +a$,
then  $c_0(H^{n}(v))  = c_0(H^{n}(u)) + a$.
The function $\Delta_n(u) = c_{0}(H^{n}(u)) - c_0(u)$ is ${\cal F}$-measurable.   We deduce that
\begin{eqnarray*}
\tilde \mu^H (H^{-n}(\{c_0 = k\}) \;\;|{\cal F} )
&= & ( \nu_{l(u)} \otimes  \mucompteurfix_{l(u)} \otimes \muoverflow)
(\{c_0 = k - \Delta_n(u)  \}  \;|\;{\cal F})\\
&= & \muunif_{2^{l_0(u)}}(k -\Delta_n(u))\\
&= &  \muunif_{2^{l_0}} (k)\\
&=& 2^{-l_0}.
\end{eqnarray*}

Remark that the last claim implies
$\mu^H (\{ c_{i}=k\} | \varphi(V_L))\ge \mu^H (H^{-n}(\{c_0 = k\}) \;\;|{\cal F} )$, so
$\mu^H (\{ c_{i}=k\} | \varphi(V_L))=2^{-l_0}$.

We prove the claim by taking the Cesaro mean and by going to the limit.
We extend the result to all integer $i$ using the shift invariance and
 the independence of $c_i$ with respect to $l_j$ when $j >i$.
We conclude by applying Lemma \ref{proprietesmesures}:
\begin{eqnarray*}
 \mu^F (\{x:c_{i}(x)=k\} |V_L )
&=&\frac{\mu^F (\{x:c_{i}(x)=k\} \cap V_L)}
{\mu^F (V_L)}\\
&=&  \frac{ \frac{L_0}{\overline{l}} \mu^H (\varphi( \{x:c_{i}(x)=k\} \cap
  V_L ))}{\frac{L_0}{\overline{l}} \mu^H (\varphi(V_L))}\\
& =&   \mu^H (\{ c_{i}=k\} | \varphi(V_L)) \\
&=& 2^{-L_i}.
\end{eqnarray*}
\hfill$\Box$
\begin{remark}
Using  Lemma 4 we obtain
$$ \muglob (\{c_{i} = k_0, \ldots , c_{i+m} = k_m\}|l_i, \ldots, l_{i+m} ) =
\left\{
\begin{array}{ll}
2^{-\sum_{j=i}^{i+m} l_j} & \hbox{ if }   \forall i\leq j \leq i+m, \, c_j < 2^{l_j},\\
0  & \hbox{otherwise.}
\end{array}
\right.
$$\end{remark}

\begin{pro}\label{shiftEntr}
The shift measurable entropy of $h_{\mu_F} (\sigma )$ is positive.
\end{pro}
\paragraph{Proof}

Let  $\alpha$  be the partition of $\Omega$ by the coordinate 0
and denote by $\alpha_{-n}^n(x)$ the element of the partition
$\alpha\vee \sigma^{-1}\alpha ...\vee \sigma^{-n+1}\vee\sigma^1
...\vee\sigma^{n}\alpha$ which contains $x$. It follows from the
definition of $\mu^F$,  that $\mu^F(\alpha_{-n}^n(x))\le \nu
\left(l(y)\in \NN^\ZZ: y\in\alpha_{-n}^n(x)\right)$. For all
$x\in\Omega$ there exist sequences of positive integers
$(-N^+_i)_{i\in\NN}$ and $(N^-_i)_{i\in\NN}$ such that
$x_{-N^+_i}=x_{N^-_i}=E$. Let
$K=\frac{\nu^3}{1-\nu}=\sum_{j>2}\nu^j$. If
 $y\in\Omega\cap\alpha_{-N^-_i}^{N^+_i}(x)$ then
 $$\nu \left(l(y)\right)=\nu_*(l_0)\pi_{k=1}^{i}\left(\nu_*(l_{-k})\nu_*(l_k)\right)
 =\nu^{N^-_i+N^+_i+1}/(K^{2i+1}).$$
 Since $l_i\ge 3$ and $\sum_{k=-i}^{i}l_k=N^-_i+N^+_i$, we obtain
$$\log \left(\nu(l(y))
\right)=(N^-_i+N^+_i+1)\log{\nu}-(2i+1)\log{K}\le
(N^-_i+N^+_i+1)\log{\frac{1}{1-\nu}}.$$ Since $\mu^F (\Omega)=1$
we get
$$
\lim_{i\to\infty}\int_{A^\ZZ} \frac{-\log (\mu^F(\alpha_{-N^-_i}^{N^+_i}(x)))}{N^-_i+N^+_i+1}d\mu^F(x)\ge
(\log{\frac{1}{1-\nu}})>0.
$$
By the probabilistic version of the Shannon-McMillan-Breiman
theorem  for a $\ZZ$-action \cite{Ki}, the left side of the
previous inequality is equal to $h_\mu^F (\sigma )$, so we can
conclude. \hfill$\Box$

\section{Lyapunov exponents}

Recall that  no information can cross a counter before it reaches
the top, that is, before time $\frac{1}{2} (\overline{c}_i - c_i)$
since at each step the counter is incremented by $1$ or $2$.
When the information reaches the next counter,
it has to wait more than $\frac{1}{2} (\overline{c}_{i+1} - c_{i+1})$,
where this quantity is estimated at the arrival time of the
information, and so on.

But each counter is uniformly distributed among its possible
values. So expectation of these times is bounded below by
$\frac{1}{4} \EE[\overline{c}_i] = \EE[2^{l_i -2}]$. A good choice
of the measure $\nu$ can make the expectation of the $l_i$ finite
so that we can define the invariant measure $\muf$. But
$\EE[2^{l_i}]$ is infinite so that expectation of time needed to
cross a counter is infinite and hence the sum of these times
divided by the sum of the length of the binary counters tends to
infinity.

Instead of being so specific, we use a rougher argument. Taking
$\mulength$ to be geometric, we show that there exists a counter
large enough to slow the speed of transmission of information. That
is, for given $n$, there is a counter of length larger than
$2\ln{n}$,  information needs a time of order $n^\delta$, with $\delta >1$ to
cross.  This is enough to conclude.

\begin{pro}\label{ii=0}
We have
$$I^+_{\muf} + I^-_{\muf}  =0.$$
\end{pro}

\paragraph{Proof}
We just have to show that $I^+_{\muf}  =0$ since $I^-_{\muf}$ is
clearly zero. Recall that $I^+_n(x)$ is the minimal number of
coordinates that we have to fix to ensure that for all $y$ such
that as soon as $y(-I_n^+(x),\infty ) = x(-I_n^+(x), \infty )$, we
have $F^k(y)(0,\infty ) = F^k(x)(0,\infty)$ for all  $0\le k\le
n$. Let
$$
t_F(n)(x) = \min \{s \; : \; \exists y, y(-n,\infty )=x(-n,\infty )
\hbox{ and } F^s(x)(0,\infty ) \neq F^s(y)(0,\infty )\}
$$
be the time needed for a perturbation to cross $n$ coordinates.
Note that $t_F(n)$ and $I^+_n$ are related by $t_F(s )(x) \geq
n\,\Leftrightarrow\,
 I^+_{n}(x)\leq s$.
We now define an analog of $t_F(n)$ for the model.
For all  $x\in \varphi (\Omega )$, let $M_n(x)$ be the smallest negative integer $m$ such that
$\sum_{i=m}^{0} l_i(x)\le n$. Define
$$
t_H(n)(x) =
\min \left\{ s \ge 0: \;  \exists y \in \Omega, \begin{array}{ll}
&y(M_n(x),\infty )=x(M_n(x),\infty ),\\
&H^s(x)(0,\infty ) \neq H^s(y)(0,\infty )
\end{array}\right\}.
$$
Define for all large enough positive integer, the subset of $\Omega$
\small
$$
U_n=\left\{x : \begin{array}{lll} &l_0(x) \leq 2\ln (n),\\
&\exists i \in \{ M_n(x), \ldots, -1\},  l_i(x) \ge 2\ln (n),\\
&c_i(x)\leq 2^{l_i(x)} - 2^{1.5\ln(n)}
\end{array} \right\}.
$$

We claim that
$$\lim_{n \to \infty} \muf (U_n) =  1.$$
Choose and fix an integer $n$ large.
Note that if $\forall i \in \{ M_n(x), \ldots, 0\}$,  $ l_i \le 2\ln (n)$,
we have $|M_n(x)|\ge \frac{n}{2\ln (n)}$. Write $M = -\lfloor
  \frac{n}{2\ln (n)}\rfloor $ and define
$$V_n =   \left\{ L \in \NN^{|M|+1} : L_0 \leq 2 \ln{n} \hbox{ and }
\forall i \in \{ M, \ldots, -1\},   L_i \le  2\ln (n) \right\}. $$
Note that on $\{l_0 \leq 2\ln{n}\}$,
existence of $i\geq M$ with $l_i \geq 2\ln{n}$ yields
existence of $i\geq M_n(x) \geq M $ with $l_i \geq 2\ln{n}$ (the same $i$).

Given a $L = (L_M, \ldots, L_0)$ in the complementary of  $V_n$ denoted by $V_n^c$, we denote
$V_L^\Omega = \left\{ x \in \Omega : (l_M(x), \ldots, l_0(x)) = L  \right\}$
and we denote $i(L)$ the larger
index $M<i<0$
with $L_{i} > 2 \ln{n}$. Let $V_n^{c*}$ be the subset of $V_n^c$ with $L_0\le 2\ln (n)$.
The measure of $U_n$ is bounded below by
\hskip -1 true cm
\small
\begin{eqnarray*}
\mu^F(U_n) & = &\sum_{L\in \NN^M} \mu^F(U_n \cap V_L^\Omega)\\
&\geq &  \sum_{L \in V_n^{c*}}
\mu^F(U_n \cap V_L^\Omega\})\\
&\geq&  \sum_{L \in V_n^{c*}}
\mu^F(\{c_{i(L)}(x) \leq 2^{l_{i(L)}} - 2^{1.5\ln{(n)}} \} \cap V_L^\Omega)\\
&=& \sum_{L \in V_n^{c*}}
\mu^F\left(c_{i(L)}(x) \leq 2^{l_{i(L)}} - 2^{1.5\ln{(n)}}  \, | \, V_L^\Omega\right)
 \; \mu^F \left( V_L^\Omega   \right).
\end{eqnarray*}
\normalsize
According to Lemma \ref{uniforme},  given the $l_i$'s, for $M
\leq i \leq 0$,  the random variable $c_i$ is uniformly
distributed. Hence for all $L \in V_n^{c*}$,
$$ \mu^F(\{x\vert c_i(x) \geq 2^{l_i} - 2^{1.5 \ln(n)}\} | V_L^\Omega)  =  2^{1.5\ln(n)}  2^{-L_i(L)}
\le 2^{1.5\ln(n) -2\ln(n)} = n^{-0.5 \ln(2)}, $$
so that
$$
\mu^F(U_n)  \geq   \sum_{L \in V_n^{c*} }
(1-n^{-0.5 \ln(2)}) \mu^F \left( V_L^\Omega   \right)
 \geq   (1-n^{-0.5\ln(2)})  \nu(V_n^{c*}). $$
Since $\muemetteur_*(l_i \ge 2\ln (n) ) =  cst. \sum_{k \geq 2 \ln(n)
} \nu^k \leq cst. n^{2 \ln{\nu}}$,
and  the $l_i$'s are independent, it is straightforward to
prove the existence of  constants $c_1$, $c_2$ and $c_3$ independent on $n$,
such that,
$$ \muemetteur \left(V_n\right) \leq  \left( 1 - \sum_{k \geq 2
    \ln(n) } \nu^k \right)^{|M|}
\leq  c_1 n^{2 \ln{\nu}} + c_2   \exp{\left( - c_3 \frac{n^{1-2\ln (\nu )}}{2\ln
     (n)}\right)}, $$
and $\nu (V_n^{c*})\geq (1-\nu (V_n))\nu_*(l_0\le 2\ln (n))$.

\paragraph{}
We conclude that
$$ \muf(U_n) \geq \left(1 - n^{-0.5\ln(2) } \right) (1-c_4n^{2 \ln{\nu}}) \left(1 -
  c_1 n^{2 \ln{\nu}} + c_2  \exp{\left( -c_3\frac{n^{1-2\ln (\nu )}}{2\ln (n)}\right)} \right).$$
Since $\nu =\frac{2}{3}$, $2\ln (\nu) <1$, this bound converges
to $0$ and the claim follows.

\paragraph{}
Now, we claim that there is a constant $c>0$ and a constant $\delta>1$
such that on $U_n$, we have
$$t_H(n)(x) \geq c  n^\delta.$$
Indeed, if $x \in U_n$, and $y$ is such that $y(-n,\infty) =
x(-n,\infty)$, then $M(x,n) = M(y,n) =: M$ and
$\varphi(y)(-M,\infty) =  \varphi(x)(-M,\infty)$. Hence  there is
an index $0 < i\leq M $ with $l_i(x) = l_i(y)= L$ and $c_i(x) =
c_i(y)= c$  satisfying
$$L \ge 2\ln (n) \hbox{ and } c \leq 2^{L} - 2^{1.5\ln(n)}.$$
Notice that $i <0$ because $x \in U_n$. For all  $s \leq
\frac{1}{2} 2^{1.5\ln(n)}$, we have $r_i^s(x) =r_i^s(y) = 0 $,
since $c_i^s (x)  < 2^{l_i(x)}$. Hence, for all $j > i$ (and in
particular for $j=0$), $c_{j}^s(x) = c_{j}^s(y)$. For $j=0$, this
implies that $t_H(n)(x) \geq \frac{1}{2} n^{1.5 \ln{2}}$. As
$1.5\ln(2)>1$ the claim holds.

\paragraph{}
There is no $y$ with $y(-n,\infty) =
x(-n,\infty)$ and $F^s(x)(0,\infty) \neq F®(y)(0,\infty)$ if
$s < t_H(n)(x)$ because  $y(-n,\infty) =
x(-n,\infty) \Rightarrow \varphi(y)(-M(x,n),\infty)
=\varphi(x)(-M(x,n),\infty)$.
This implies  that $t_F(n)(x) \geq t_H(n)(x)$.
It follows that $t_F(n)(x) \geq c n^\delta$  on $U_n$.

\paragraph{}
Setting $s = \left[ n^{\frac{1}{\delta}}\right]$, we see that
$ t_F(\left[ n^{\frac{1}{\delta}}\right]) \geq n \Leftrightarrow I^+_{n} \leq
\left[ n^{\frac{1}{\delta}}\right].$
We deduce that there is a constant $c$ such that, on $U_n$,
$$\frac{I^+_n(x)}{n} \leq c n^{\frac{1}{\delta} -1}. $$
Since $\frac{I^+_n(x)}{n}$ is bounded (by $r=2$, radius of the
automaton),
the conclusion follows
from the inequality
$$\int \frac{I^+_n}{n}d\mu^F \leq \int_{U_n} c \,
 n^{\frac{1}{\delta} -1}    d \mu^F  +  2 \mu^F(U_n^c) .$$
That is,
$$I^+_{\mu^F}=0.$$
\hfill$\Box$

Using Proposition \ref{i=0} we obtain
\begin{cor}
The measurable entropy $h_{\mu^F} (F)$ is equal to zero.
\end{cor}

\section{Questions}
We end this paper with a few open questions and conjectures.

\begin{conj}
The measure $\mu_F$  is shift-ergodic.
\end{conj}
The measure $\mu^F$ is clearly not $F$-ergodic since the Es do not
move. It is still not clear to us  whether or not it is possible to
construct a sensitive automaton with null Lyapunov exponents for a
$F$-ergodic measure.

\begin{conj}
A sensitive cellular automaton acting surjectively on an irreducible subshift of finite type has average
positive Lyapunov exponents if the invariant measure we consider is the Parry measure on this subshift.
\end{conj}
\begin{conj}
If a cellular automaton has no equicontinuous points (i.e. it is
sensitive), then there exists a point $x$ such that $\liminf
\frac{I^+_n(x)}{n}>0$ or $\liminf
\frac{I^-_n(x)}{n}>0$.
\end{conj}



\section{References}


\end{document}